\documentclass[11pt]{amsart}

\usepackage{amssymb,amsthm,amsmath}
\usepackage[numbers,sort&compress]{natbib}
\usepackage{color}
\usepackage{graphicx}
\usepackage{tikz}
\usepackage[colorlinks,
linkcolor=red,
anchorcolor=blue,
citecolor=green
]{hyperref}

\def\a{\alpha}

\let\newpf\proof \let\proof\relax 
\newenvironment{pf}{\newpf[\proofname]}{\qed\endtrivlist}

\newcommand{\ba}{\overline{A}}

\def\be{\begin{equation}}
\def\ee{\end{equation}}

\def\ba{{\begin{align}}}
\def\ea{{\end{align}}}

\def\bm{\begin{matrix}}
\def\em{\end{matrix}}

\def\u{{\mathbb U}}

\def\a{{\alpha}}

\def\0{{\mathbf 0}}

\newtheorem{Theorem}{Theorem}[section]
\newtheorem{Lemma}{Lemma}[section]
\newtheorem{Proposition}{Proposition}[section]

\newtheorem{Remark}{Remark}[section]

\numberwithin{equation}{section}

\theoremstyle{definition}

\newcommand{\C}{{\mathbb C}}

\newcommand{\N}{{\mathbb N}}

\newcommand{\R}{{\mathbb R}}
\newcommand{\T}{{\mathbb T}}

\newcommand{\Z}{{\mathbb Z}}

\def\B0{{\bold{0}}}

\catcode`\@=12

\def\Empty{}
\newcommand\oplabel[1]{
  \def\OpArg{#1} \ifx \OpArg\Empty {} \else
    \label{#1}
  \fi}

\newcommand{\comm}[1]{}
\newcommand{\comment}[1]{}

\begin{document}

\title[]{H\"older Continuity of Absolutely Continuous Spectral Measure for Multi-frequency Schr\"odinger Operators}
\author{Xin Zhao}
\address{
Department of Mathematics, Nanjing University, Nanjing 210093, China
}

 \email{njuzhaox@126.com}

\begin{abstract}
We establish sharp results on the modulus of continuity of the distribution of the spectral measure for multi-frequency Schr\"odinger operators with Diphantine frequencies and small analytic potentials.
\end{abstract}
\maketitle

\section{\textbf{Introduction}}
In this paper we consider multi-frequency quasi-periodic
Schr\"odinger operator on $\ell^2(\Z)$,
\begin{equation}\label{1}
(H_{\lambda V,\alpha,\theta}u)_n=u_{n+1}+u_{n-1}+ \lambda V(\theta+n\alpha)u_n,\ \ n\in\Z,
\end{equation}
where $\alpha,\theta\in \T^d$ are parameters (called the frequency, phase respectively), $V\in C^\omega(\T^d,\R)$ is called the potential and $\lambda\in\R$ is called the coupling constant.

There are two fundamental quantities in the study of the spectral theory of quasi-periodic Schr\"odinger operators, the {\it Lyapunov exponent} (LE) and the {\it integrated density of states} (IDS)\footnote{One can consult Section 2.3 for the definitions.}. On the one hand, it is well understood that the regularity of LE (IDS) plays an important role in the study of spectral theory of quasi-periodic Schr\"odinger operators. For example, the classical Kotani theory \cite{kotani} says that if the Lyapunov exponent vanishes in the spectrum, then the absolutely continuity of IDS is equivalent to that the absolute continuity of the spectral measure for {\it a.e.} $\theta$. Recently, it is explored in \cite{dgl,dgsv} the H\"older continuity of LE (IDS) plays an important role in certain topological structure (called {\it Homogeneity}) of spectrum set of quasi-periodic Schr\"odinger operators. On the other hand, the regularity of LE itself is one of the fundamental questions in dynamical systems, for partial results, one can see Viana \cite{viana} and references therein. The Lyapunov exponent is connected to the integrated density of states by the {\it Thouless formula} which was derived on a non-rigorous basis by Thouless \cite{th} and then rigorously proved by Avron and Simon \cite{as}. Classical hard analysis \cite{gs1} indicates that this formula transfers the H\"older regularity of IDS to that of LE. Thus the regularity of LE (IDS) reduce to a problem of regularity of IDS.

For quasi-periodic Schr\"odinger operators, let $\mu_\theta=\mu^{\delta_0}_{\lambda V,\alpha,\theta}+\mu^{\delta_1}_{\lambda V,\alpha,\theta}$ be the associated universal spectral measure (see Section 2.1 for details), it is standard that the IDS is the average of $\mu_\theta$ in $\theta$. Thus, a more difficult and subtle question is the regularity of the distribution of individual spectral measure $\mu_\theta$. Note that if $\mu_\theta=\mu_\theta^{pp}$, then the distribution of $\mu_\theta$ is not even continuous. Thus a more suitable question is the regularity of the distribution of $\mu_{\theta}$ when $\mu_{\theta}=\mu_{\theta}^{ac}$. This question was answered by Avila and Jitomirskaya \cite{aj} in the one-frequency case under the assumption that the frequency is Diophantine \footnote{ $\alpha \in\T^d$ is  called {\it Diophantine}, denoted by $\alpha \in {\rm DC}_d(\kappa,\tau)$, if there exist $\kappa>0$ and $\tau>d-1$ such that
\begin{equation}\label{dio}
{\rm DC}_d(\kappa,\tau):=\left\{\alpha \in\T^d:  \inf_{j \in \Z}\left| \langle n,\alpha  \rangle - j \right|
> \frac{\kappa}{|n|^{\tau}},\quad \forall \  n\in\Z^d\backslash\{0\} \right\}.
\end{equation}
Let ${\rm DC}_d:=\bigcup_{\kappa>0,
\, \tau>d-1} {\rm DC}_d(\kappa,\tau)$.}. And it was generalized by Liu and Yuan \cite{ly} to Liouvillean frequency.

However, the results in \cite{aj,ly} were restricted to one-frequency case, since a crucial technique in \cite{aj,ly} is {\it almost reducibility} developed by Avila and Jitomirskaya in \cite{avilajitomirskaya} based on quantitative Aubry duality and it seems non-trivial to generalize the method in \cite{avilajitomirskaya} to multi-frequency case. While almost reducibility can also be got directly by classical KAM theory \cite{ds,e1,hy,lyzz,ccyz}, and results in \cite{ds,e1,lyzz,ccyz} do work in any dimension. This allows the possibility to study the regularity of distribution of absolutely continuous spectral measure in the multi-frequency case.

In this paper, we generalize the results in \cite{aj} to multi-frequency case by the KAM scheme recently developed in \cite{lyzz,ccyz} for multi-frequencies. Our main Theorem is
\begin{Theorem}\label{main}
Assume $\alpha\in DC_d$ and $V\in C^\omega(\T^d,\R)$, there exist $\lambda_0(\alpha,V)$ and $C=C(\alpha,\lambda V)$ such that if $\lambda<\lambda_0$, then for any $E\in \R$, there holds
$$
\mu_{\theta}(E-\epsilon,E+\epsilon)\leq C\epsilon^{\frac{1}{2}},
$$
for any $\epsilon>0$ and any $\theta\in\T^d$.
\end{Theorem}
\begin{Remark}
The $\frac{1}{2}$-modulus of continuity is sharp since even for IDS (the average), this is sharp \cite{p}.
\end{Remark}

We now give a brief review of the histories on the regularity of LE (IDS) for quasi-periodic Schr\"odinger operators. For analytic quasi-periodic potentials, in the positive Lyapunov exponent regime, Goldstein and Schlag \cite{gs1} developed some sharp version of large deviation theorems for real analytic potentials with strong Diophantine frequency, moreover they further developed Avalanche Principle and proved that LE is H\"{o}lder continuous (one-frequency) or weak H\"{o}lder continuous (multi-frequency). For the Almost Mathieu operator where $V(\theta)=\lambda \cos2\pi(\theta)$ and $d=1$, Bourgain \cite{bourgain} proved that for Diophantine $\alpha$ and large enough $\lambda$, LE is $\frac{1}{2}-\epsilon$-H\"{o}lder continuous for any $\epsilon>0$.  Later, Goldstein and Schlag \cite{goldsteinschlag2} generalized Bourgain's result \cite{bourgain}, and proved that if the potential  is in a small $L^{\infty}$ neighborhood of a trigonometric polynomial of degree $k$, then the IDS is H\"older $\frac{1}{2k}-\epsilon$-continuous for all $\epsilon>0$. Moreover, they  further proved (\cite{goldsteinschlag2} ) that IDS is absolutely continuous for $a.e.$ $\alpha$.

In the zero Lyapunov exponent regime, based on Eliasson's perturbative KAM scheme \cite{e1}, Amor \cite{amor} got $\frac{1}{2}$-H\"{o}lder continuity of IDS for quasi-periodic cocycles in $SL(2,\R)$ with Diophantine frequency. Besides, Avila and Jitomirskaya \cite{avilajitomirskaya} used almost localization and Aubry duality to obtain the same result with one frequency in the non-perturbative regime. It is worth mentioning that all the above results require Diophantine or strong Diophantine conditions. For small potentials and generic frequencies, it is actually possible to show that the Lyapunov exponent is not  H\"{o}lder continuous. A recent breakthrough belongs to Avila \cite{avila}: for one-frequency Schr\"{o}dinger operators with general analytic potentials and irrational frequency, Avila \cite{avila} has established the fantastic global theory saying that  Lyapunov exponent is a $C^\omega$-stratified function of the energy.

For the lower regularity case, Klein \cite{klein} proved that for Schr\"{o}dinger operators with potentials in a Gevrey class, the Lyapunov exponent is weak H\"{o}lder continuous on any compact interval of the energy provided that the coupling constant is large enough, the frequency is Diophantine and the potential satisfies some transversality condition.
Recently, Wang and Zhang \cite{wangzhang} obtained the weak H\"{o}lder continuity of Lyapunov exponent as a function of energies, for a class of $C^2$ quasi-periodic potentials and for any Diophantine frequency. More recently, Cai, Chavaudret, You and Zhou \cite{ccyz} proved sharp H\"older continuity of Lyapunov for quasi-periodic Schr\"odinger operator with small finitely differential potentials and Diophantine frequencies.

Up to now, the regularity result of the distribution of individual spectral measure for quasi-periodic Schr\"odinger operators is few. We mention that recently, Avila and Jitomirskaya \cite{aj} proved sharp H\"older continuity of $\mu_\theta$ in the non-perturbative regime for Diophantine frequencies. Later, Liu and Yuan generalized this result to Liouvillean frequencies.

Finally, we give the structure of this paper. Several preparation propositions and basic concepts are given in Section 2. In Section 3, we derive some quantitative estimates based on the KAM scheme developed in \cite{lyzz,ccyz}. In Section 4, based on these quantitative estimates on conjugation transformation and constant matrix, we will give the proof of the main theorem with some arguments in \cite{aj}.
\section{\textbf{Preliminary}}

We denote $C_r^\omega({\T}^d,*)$  by the space of analytic matrix-valued functions with analytic radius $r>0$. Here $``*"$ can be $\R,\C,SL(2,\R),sl(2,\R),SL(2,\C)$ and $sl(2,\C)$. The norms are defined as $$\lVert F \rVert _{r}=\sup_{\substack{
                             |\Im \theta|< r
                          }}\lVert F(\theta) \rVert.
$$
We denote by
$$
C^\omega({\T}^d,*)=\cup_{r>0}C_r^\omega({\T}^d,*),
$$
$$
\|F\|_0=\sup_{\theta\in\T^d}\|F(\theta)\|.
$$
\subsection{\textbf{Spectral measure for one dimension Schr\"odinger operator}}
\noindent

Given a bounded map $V:\Z\rightarrow \R$ called the potential, we define the associated Schr\"odinger operator by
\begin{equation}\label{schro}
(Hu)_n=u_{n+1}+u_{n-1}+V(n)u_n.
\end{equation}
It is easy to check that $H$ is a bounded self-adjoint operator in $\ell^2(\Z)$.

It is standard that for any compactly supported $\phi\in\ell^2(\Z)$ can be written as $\phi=p(H)\delta_1+q(H)\delta_0$ with suitable polynomials $p(\cdot)$ and $q(\cdot)$. Using this observation, it can then be shown that
$$
\mu=\mu_{\delta_0}+\mu_{\delta_1}
$$
where $\mu_{\delta_0}$, $\mu_{\delta_1}$ are the associated spectral measures of $H$ with respect to $\delta_0$, $\delta_1$, can serve as a universal spectral measure for $H$. More precisely, for any $\phi\in\ell^2(\Z)$, $\mu_\phi$ is absolutely continuous with respect to $\mu$.

\subsection{Borel transformation of spectral measure and $m$ function}
\noindent

The Borel transform of $\mu$ takes the form
\begin{equation}\label{M}
M(z):=F_\mu(z)=\int\frac{d\mu(E)}{E-z}=\langle\delta_0,(H-zI)^{-1}\delta_0\rangle+\langle\delta_1,(H-zI)^{-1}\delta_1\rangle.
\end{equation}
It is standard that $M(z)$ has a close relation to the well-known Weyl-Titchmarsh $m$-function. Given $z\in\C^+$, then there are non-zero solutions $u_z^{\pm}$ of $Hu_z^{\pm}=zu_z^{\pm}$ which are $\ell^2$ at $\pm\infty$. The Weyl-Titchmarsh $m$-functions are defined by
$$
m_z^{\pm}=\mp\frac{u_z^{\pm}(1)}{u_z^{\pm}(0)},
$$
we refer readers to consult \cite{jl2,jl3} for more details of the definition.

As discussed in \cite{jl3},
$$
M(z)=\frac{m_z^+m_z^--1}{m_z^++m_z^-}.
$$

For $k\in\Z^+$ and $E\in\R$, let
$$
P_k(E)=\sum\limits_{j=1}^k A^*_{2j-1}(E)A_{2j-1}(E),
$$
where $A_n(E)=T(E,n)\cdots T(E,1)T(E,0)$ with $T(E,n)=\begin{pmatrix}E-V(n)&-1\\1&0 \end{pmatrix}$.

 The following propositions which were proved in \cite{aj} are very important for our applications. For completeness, we give the proof here.
\begin{Proposition}\label{mfunction}
For any $E\in\R$ and $\epsilon_k=\sqrt{\frac{1}{4\det{P_k(E)}}}$, we have
$$
\mu(E-\epsilon_k,E+\epsilon_k)\leq 2\epsilon_k\Im M(E+i\epsilon_k)\leq 4(5+\sqrt{24})\epsilon_k^2\|P_k(E)\|.
$$
\end{Proposition}
\begin{pf}
By \eqref{M}, for any $\epsilon>0$, we have
\begin{align}\label{Mfunction}
\Im M(E+i\epsilon)=\int\frac{\epsilon}{(E'-E)^2+\epsilon^2}d\mu(E'),
\end{align}
thus
\begin{align}\label{ine1}
\Im M(E+i\epsilon)&\geq \int_{E-\epsilon}^{E+\epsilon}\frac{\epsilon}{(E'-E)^2+\epsilon^2}d\mu(E')\\ \nonumber
&\geq \frac{1}{2\epsilon}\int_{E-\epsilon}^{E+\epsilon}d\mu(E')\\ \nonumber
&=\frac{1}{2\epsilon}\mu(E-\epsilon,E+\epsilon).
\end{align}
We denote by $\psi(z)=\sup_\beta|R_{-\beta/2\pi}z|$, by the argument in Section 4.1 in \cite{aj}, one has
\begin{equation}\label{ine2}
|M(E+i\epsilon)|\leq \psi(m_{E+i\epsilon}^+).
\end{equation}
By Lemma 4.2 in \cite{aj}, one has
\begin{equation}\label{ine3}
\psi(m_{E+i\epsilon_k}^+)\leq 2(5+\sqrt{24})\epsilon_k\|P_k(E)\|.
\end{equation}
\eqref{ine1}, \eqref{ine2} and \eqref{ine3} imply that
\begin{align*}
\mu(E-\epsilon_k,E+\epsilon_k)&\leq 2\epsilon_k\Im M(E+i\epsilon_k)\leq 2\epsilon_k\psi(m_{E+i\epsilon_k}^+)\\
&\leq 4(5+\sqrt{24})\epsilon_k^2\|P_k(E)\|.
\end{align*}
\end{pf}
\begin{Proposition}\label{prob}
For any $E\in\R$ and $\epsilon_k=\sqrt{\frac{1}{4\det{P_k(E)}}}$, let $(u_n^\beta)_{n\geq 0}$ satisfy
$$
T(E,n)\begin{pmatrix}u_n^\beta\\u_n^\beta\end{pmatrix}=\begin{pmatrix}u_{n+1}^\beta\\u_{n+1}^\beta\end{pmatrix},
$$
$$
u_0^\beta\cos\beta+u_1^\beta\sin\beta=0,\ \  |u_0^\beta|^2+|u_1^\beta|^2=1.
$$
Then we have
$$
\det{P_k(E)}=\inf\limits_{\beta}\|u^\beta\|^2_{2k}\|u^{\beta+\pi/2}\|^2_{2k},
$$
where for any integer $L$
$$
\|u\|_{L}=(\sum\limits_{n=1}^L|u_n|^2)^{\frac{1}{2}}.
$$
\end{Proposition}
\begin{pf}
By the definition of $P_k(E)$, we have
$$
\|u^\beta\|^2_{2k}=\langle P_k(E)\begin{pmatrix}u_1^\beta\\u_0^\beta\end{pmatrix},\begin{pmatrix}u_1^\beta\\u_0^\beta\end{pmatrix}\rangle.
$$
Since $P_k(E)$ is self-adjoint, it immediately follows that
$$
\det{P_k(E)}=\inf\limits_{\beta}\|u^\beta\|^2_{2k}\|u^{\beta+\pi/2}\|^2_{2k}.
$$
\end{pf}
\subsection{Linear algebra preparations}
\noindent

In this subsection, we give some basic facts in linear algebra. Note that, essentially, all these facts are proved in \cite{aj}. We give the proof here for completeness.
\begin{Proposition}\label{comp}
Assume that
$$
T=\begin{pmatrix} e^{2\pi i\theta}&c\\0&e^{-2\pi i\theta}
\end{pmatrix},
$$
let $X_k=\sum\limits_{j=1}^k (T^{2j-1})^*(T^{2j-1})$, then
$$
X_k=\begin{pmatrix}
k&x_{k,1}\\
\bar{x}_{k,1}&x_{k,2}
\end{pmatrix}
$$
where
$$
x_{k,1}=ce^{-2\pi i\theta}\sum\limits_{j=1}^k\frac{e^{-4\pi i\theta(2j-1)}-1}{e^{-4\pi i \theta}-1},
$$
$$
x_{k,2}=k+|c|^2\sum\limits_{j=1}^k(\frac{\sin2\pi(2j-1)\theta}{\sin2\pi\theta})^2.
$$
\end{Proposition}
\begin{pf}
We prove this by induction. If $k=1$, direct computation shows that
\begin{align*}
X_1&=\begin{pmatrix} e^{-2\pi i\theta}&0\\\bar{c}&e^{2\pi i\theta}
\end{pmatrix}\begin{pmatrix} e^{2\pi i\theta}&c\\0&e^{-2\pi i\theta}
\end{pmatrix}\\
&=\begin{pmatrix} 1&ce^{-2\pi i\theta}\\ \bar{c}e^{2\pi i\theta}&1+|c|^2
\end{pmatrix}.
\end{align*}
Assume that we are at the $n$-th step, we consider the $(n+1)$-th step, note that
$$
T^{2n+1}=\begin{pmatrix}
e^{2\pi i(2n+1)\theta}&t_{2n+1}\\
0&e^{-2\pi i(2n+1)\theta},
\end{pmatrix}
$$
where $t_{2n+1}=ce^{4\pi in\theta}\frac{e^{-4\pi i(2n+1)\theta}-1}{e^{-4\pi i\theta}-1}$. Thus we have
\begin{align*}
X_{n+1}&=X_n+(T^{2n+1})^*T^{2n+1}\\
&=\begin{pmatrix}
n&x_{n,1}\\
\bar{x}_{n,1}&x_{n,2}
\end{pmatrix}
+\begin{pmatrix} 1&t_{2n+1}e^{-2\pi i(2n+1)\theta}\\ \bar{t}_{2n+1}e^{2\pi i(2n+1)\theta}&1+|t_{2n+1}|^2
\end{pmatrix},
\end{align*}
this implies that
\begin{align*}
x_{n+1,1}&=x_{n,1}+t_{2n+1}e^{-2\pi i(2n+1)\theta}\\
&=ce^{-2\pi i\theta}\sum\limits_{j=1}^n\frac{e^{-4\pi i\theta(2j-1)}-1}{e^{-4\pi i \theta}-1}+ce^{4\pi in\theta}\frac{e^{-4\pi i(2n+1)\theta}-1}{e^{-4\pi i\theta}-1}e^{-2\pi i(2n+1)\theta}\\
&=ce^{-2\pi i\theta}\sum\limits_{j=1}^{n+1}\frac{e^{-4\pi i\theta(2j-1)}-1}{e^{-4\pi i \theta}-1}.
\end{align*}
\begin{align*}
x_{n+1,2}&=x_{n,2}+1+|t_{2n+1}|^2\\
&=n+|c|^2\sum\limits_{j=1}^n(\frac{\sin2\pi(2j-1)\theta}{\sin2\pi\theta})^2+1+|c|^2\frac{|e^{-4\pi i(2n+1)\theta}-1|^2}{|e^{-4\pi i\theta}-1|^2}\\
&=n+1+|c|^2\sum\limits_{j=1}^{n+1}(\frac{\sin2\pi(2j-1)\theta}{\sin2\pi\theta})^2.
\end{align*}
\end{pf}
\subsection{\textbf{Cocycles, Lyapunov exponents  and fibered rotation number}}
\noindent

Given $A \in C^0(\T^d,{\rm SL}(2,\C))$ and rationally independent $\alpha\in\R^d$, we define the quasi-periodic \textit{cocycle} $(\alpha,A)$:
$$
(\alpha,A)\colon \left\{
\begin{array}{rcl}
\T^d \times \C^2 &\to& \T^d \times \C^2\\[1mm]
(x,v) &\mapsto& (x+\alpha,A(x)\cdot v)
\end{array}
\right.  .
$$
The iterates of $(\alpha,A)$ are of the form $(\alpha,A)^n=(n\alpha,  \mathcal{A}_n)$, where
$$
\mathcal{A}_n(x):=
\left\{\begin{array}{l l}
A(x+(n-1)\alpha) \cdots A(x+\alpha) A(x),  & n\geq 0\\[1mm]
A^{-1}(x+n\alpha) A^{-1}(x+(n+1)\alpha) \cdots A^{-1}(x-\alpha), & n <0
\end{array}\right.    .
$$
The {\it Lyapunov exponent} is defined by
$\displaystyle
L(\alpha,A):=\lim\limits_{n\to \infty} \frac{1}{n} \int_{\T^d} \ln \|\mathcal{A}_n(x)\| dx
$.

The cocycle $(\alpha,A)$ is {\it uniformly hyperbolic} if, for every $x \in \T^d$, there exists a continuous splitting $\C^2=E^s(x)\oplus E^u(x)$ such that for every $n \geq 0$,
$$
\begin{array}{rl}
|\mathcal{A}_n(x) \, v| \leq C e^{-cn}|v|, &  v \in E^s(x),\\[1mm]
|\mathcal{A}_n(x)^{-1}   v| \leq C e^{-cn}|v|, &  v \in E^u(x+n\alpha),
\end{array}
$$
for some constants $C,c>0$.
This splitting is invariant by the dynamics, i.e.,
$$A(x) E^{*}(x)=E^{*}(x+\alpha), \quad *=``s"\;\ {\rm or} \;\ ``u", \quad \forall \  x \in \T^d.$$

Assume that $A \in C^0(\T^d, {\rm SL}(2, \R))$ is homotopic to the identity. It induces the projective skew-product $F_A\colon \T^d \times \mathbb{S}^1 \to \T^d \times \mathbb{S}^1$ with
$$
F_A(x,w):=\left(x+\a,\, \frac{A(x) \cdot w}{|A(x) \cdot w|}\right),
$$
which is also homotopic to the identity.
Thus we can lift $F_A$ to a map $F_A\colon \T^d \times \R \to \T^d \times \R$ of the form $F_A(x,y)=(x+\alpha,y+\psi_x(y))$, where for every $x \in \T^d$, $\psi_x$ is $\Z$-periodic.
The map $\psi\colon\T^d \times \T  \to \R$ is called a {\it lift} of $A$. Let $\mu$ be any probability measure on $\T^d \times \R$ which is invariant by $F_A$, and whose projection on the first coordinate is given by Lebesgue measure.
The number
$$
\rho(\alpha,A):=\int_{\T^d \times \R} \psi_x(y)\ d\mu(x,y) \ {\rm mod} \ \Z
$$
 depends  neither on the lift $\psi$ nor on the measure $\mu$, and is called the \textit{fibered rotation number} of $(\alpha,A)$ (see \cite{H,JM} for more details).

Given $\theta\in\T^d$, let $
R_\theta:=
\begin{pmatrix}
\cos2 \pi\theta & -\sin2\pi\theta\\
\sin2\pi\theta & \cos2\pi\theta
\end{pmatrix}$.
If $A\colon \T^d\to{\rm PSL}(2,\R)$ is homotopic to $\theta \mapsto R_{\frac{\langle n, \theta\rangle}{2}}$ for some $n\in\Z^d$,
then we call $n$ the {\it degree} of $A$ and denote it by $\deg A$.
The fibered rotation number is invariant under real conjugacies which are homotopic to the identity. More generally, if $(\alpha,A_1)$ is conjugated to $(\alpha, A_2)$, i.e., $B(\cdot+\alpha)^{-1}A_1(\cdot)B(\cdot)=A_2(\cdot)$, for some $B \colon \T^d\to{\rm PSL}(2,\R)$ with ${\rm deg} B=n$, then
\begin{equation}\label{rotation number}
\rho(\alpha, A_1)= \rho(\alpha, A_2)+ \frac{\langle n,\alpha \rangle}2.
\end{equation}

%

A typical  example is given by the so-called \textit{Schr\"{o}dinger cocycles} $(\alpha,S_E^{V})$, with
$$
S_E^{V}(\cdot):=
\begin{pmatrix}
E-V(\cdot) & -1\\
1 & 0
\end{pmatrix},   \quad E\in\R.
$$
Those cocycles were introduced because it is equivalent to the eigenvalue equation $H_{V, \alpha, \theta}u=E u$. Indeed, any formal solution $u=(u_n)_{n \in \Z}$ of $H_{V, \alpha, \theta}u=E u$ satisfies
$$
\begin{pmatrix}
u_{n+1}\\
u_n
\end{pmatrix}
= S_E^V(\theta+n\alpha) \begin{pmatrix}
u_{n}\\
u_{n-1}
\end{pmatrix},\quad \forall \  n \in \Z.
$$
The spectral properties of $H_{V,\alpha,\theta}$ and the dynamics of $(\alpha,S_E^V)$ are closely related by the well-known fact:
 $E\in \Sigma_{\alpha,V}$ if and only if $(\alpha,S_E^{V})$ is \textit{not} uniformly hyperbolic. Throughout the paper, we will denote $L(E)=L(\alpha,S_E^{V})$  and $\rho(E)=\rho(\alpha,S_E^{V})$ for short.

It is well known that the spectrum of $H_{V,\alpha,\theta}$ denote by $\Sigma_{\alpha,V}$, is a compact subset of $\R$, independent of $\theta$ if $(1,\alpha)$ is rationally independent. The {\it integrated density of states} (IDS) $N_{\alpha,V}:\R\rightarrow [0,1]$ of $H_{V,\alpha,\theta}$ is defined as
$$
N_{\alpha,V}(E):=\int_{\T}\mu_{V,\alpha,\theta}(-\infty,E]d\theta,
$$
where $\mu_{V,\alpha,\theta}$ is the spectral measure of $H_{V,\alpha,\theta}$.

It is also known that $\rho(E)\in[0,\frac{1}{2}]$ relates to the integrated density of states $N=N_{\alpha,V}$ as follows:
$$
N(E)=1-2\rho(E).
$$

\section{Quantitative almost reducibility}
In this section, we concentrate on the following analytic quasi-periodic $SL(2,\R)$ cocycle:
$$
(\alpha,A_0e^{f_0(\theta)}):\T^{d}\times\R^{2} \rightarrow \T^{d}\times\R^{2};(\theta,v)\mapsto (\theta+\alpha,A_0e^{f_0(\theta)}\cdot v),
$$
where
$f_0\in C^{\omega}_{r_0}(\T^{d},sl(2,\R))$,  $r_0>0$, $d\in \N^+$£¬
and $\alpha\in DC_d$.
Notice that A has eigenvalues $\{e^{i\xi},e^{-i\xi}\}$ with $\xi\in \C$.

We will prove the following quantitative almost reducibility proposition.
\begin{Proposition}\label{reducibility}
For any $0<r<r_0$, $\kappa>0$, $\tau>d-1$.
Suppose that $\alpha\in DC_d(\kappa,\tau)$. Then there exist $B_n\in C_{r}^\omega(\T^d, PSL(2,\R))$ and $A_n\in SL(2,\R)$ satisfying
$$
B_n^{-1}(\theta+\alpha)A_0e^{f_0(\theta)}B_n(\theta)=A_ne^{f_n(\theta)},
$$
provided that $\|f_0\|_{r_0}<\epsilon_*$ for some $\epsilon_*>0$ depending on $A_0,\kappa,\tau,r,r_0,d$, with the following estimates
\begin{equation}\label{est1}
\|f_n\|_{r}\leq \epsilon_n,
\end{equation}
\begin{equation}\label{est2}
\|B_n\|_{0}\leq \epsilon_{n-1}^{-\frac{1}{800}}.
\end{equation}
Moreover, there exists unitary $U_n\in SL(2,\C)$ such that
$$
U_nA_{n}U_n^{-1}=\begin{pmatrix} e^{i\xi_n} & c_n\\ 0 & e^{-i\xi_n} \end{pmatrix},
$$
and
\begin{equation}\label{estisharp}
\lvert c_n\rvert \lVert B_{n}\rVert_0^8\leq 4\lVert A_0\rVert,
\end{equation}
with $\xi_n,c_n\in \C$.
\end{Proposition}
\begin{pf}
 We prove Proposition \ref{reducibility} by iteration. Suppose that
\begin{align*}
\|f_0\|_{r_0}\leq \epsilon_* \leq \frac{c}{\|A_0\|^{D}}(r_0-r)^{D\tau},\end{align*}
where $c$, $D$ are defined in Proposition \ref{prop1}. Then we can define the sequence inductively.  Let
$\epsilon_0=\epsilon_*$, assume that we are at the $(j+1)^{th}$ KAM step, i.e. we already construct $B_j\in C^\omega_{r_{j}}(\T^d,PSL(2,\R))$ such that
$$
B_{j}^{-1}(\theta+\alpha)A_0e^{f_0(\theta)}B_{j}(\theta)=A_{j}e^{f_{j}(\theta)},
$$
where $A_j\in SL(2,\R)$ with two eigenvalues $e^{\pm i\xi_j}$ and
\begin{equation*}
\|f_j\|_{r_j}\leq \epsilon_j\leq \epsilon_0^{2^j},\ \ \|B_j\|_{0}\leq \epsilon_{j-1}^{-\frac{1}{800}}.
\end{equation*}
Then we define
$$
r_j-r_{j+1}=\frac{r_0-r}{4^{j+1}}, \ \ N_j=\frac{2|\ln\epsilon_j|}{r_j-r_{j+1}}.
$$

By our selection of  $\epsilon_0$, one can check that
\begin{equation}\label{iter}
\epsilon_j \leq \frac{c}{\|A_j\|^{D}}(r_j-r_{j+1})^{D\tau}.
\end{equation}
 Indeed, $\epsilon_j$ on the left side of the inequality decays at least super-exponentially with $j$, while $(r_j-r_{j+1})^{D\tau}$ on the right side decays exponentially with $j$.

Note that $(\ref{iter})$ implies that  Proposition \ref{prop1} can be applied iteratively, consequently one can construct
$$
\bar{B}_j\in C^\omega_{r_{j+1}}(\T^d,PSL(2,\R)),\ \ A_{j+1}\in SL(2,\R),\ \ f_{j+1}\in C_{r_{j+1}}^{\omega}(\T^d,sl(2,\R))
$$
such that
$$
\bar{B}_j^{-1}(\theta+\alpha)A_je^{f_j(\theta)}\bar{B}_j(\theta)=A_{j+1}e^{f_{j+1}(\theta)}.
$$
More precisely, we can distinguish two cases:\\

\noindent \textbf{Non-resonant case:}  If for any $n\in \Z^{d}$ with $0< |n| \leq N_j$, we have
$$
\| 2\xi_j - <n,\alpha> \|_{\R/\Z}\geq \epsilon_j^{\frac{1}{10}},
$$
then
\begin{equation*}
\| \bar{B}_j-Id\|_{r_{j+1}}\leq \epsilon_j^{\frac{1}{2}} ,\   \ \| f_{j+1}\|_{r_{j+1}}\leq \epsilon_j^2:=  \epsilon_{j+1}, \ \ \|A_{j+1}-A_j\|\leq 2\|A_j\|\epsilon_j.
\end{equation*}
Let $B_{j+1}=B_j(\theta)\bar{B}_j(\theta)$, we have
$$
B_{j+1}^{-1}(\theta+\alpha)A_0e^{f_0(\theta)}B_{j+1}(\theta)=A_{j+1}e^{f_{j+1}(\theta)},
$$
with estimate $$
\|B_{j+1}\|_{0}\leq 2\|B_j\|_0\leq 2\epsilon_{j-1}^{-\frac{1}{800}}\leq \epsilon_j^{-\frac{1}{800}}.
$$

%

\noindent \textbf{Resonant case:} If there exists $n_j$ \footnote{We call such $n_j$ the resonance.} with $0<| n_j| \leq N_j$ such that
$$
\| 2\xi_j- <n_j,\alpha> \|_{\R/\Z}< \epsilon_j^{\frac{1}{10}},
$$
then
$$\|\bar{B}_j\|_{0}\leq\epsilon_{j}^{-\frac{1}{1600}},\ \  \| f_{j+1}\|_{r_{j+1}} \ll \epsilon_j^{1600}:= \epsilon_{j+1}.$$
Moreover, $A_{j+1}=e^{A_{j+1}''}$ with $\lVert A_{j+1}''\rVert \leq 2\epsilon_j^{\frac{1}{10}}$.

Let $B_{j+1}(\theta)=B_j(\theta)\bar{B}_j(\theta)$, then we have
$$
B_{j+1}^{-1}(\theta+\alpha)A_0e^{f_0(\theta)}B_{j+1}(\theta)=A_{j+1}e^{f_{j+1}(\theta)},
$$
with
\begin{align*}
\|B_{j+1}\|_{0}&\leq \epsilon_j^{-\frac{1}{1600}}\epsilon_{j-1}^{-\frac{1}{800}}\leq \epsilon_j^{-\frac{1}{800}}.
\end{align*}

Finally, we give the proof of \eqref{estisharp}, If the $n$-th step is in the resonant case, we have
$$
A_{n}=e^{A_{n}''}, \ \ \lVert A_{n}''\rVert < 2\epsilon_{n-1}^{\frac{1}{10}}.
$$
Thus
$$
\lVert A_{n}\rVert\leq 1+4\epsilon_{n-1}^{\frac{1}{10}}\leq 2\lVert A_0\rVert,
$$
then there exists unitary $U_n\in SL(2,\C)$ such that
\begin{equation}
\label{estt}U_nA_{n}U_n^{-1}=\begin{pmatrix} e^{i\xi_{n}} & c_{n}\\ 0 & e^{-i\xi_{n}} \end{pmatrix},
\end{equation}
with $\lvert c_{n}\rvert\leq 2\lVert A_{n}''\rVert \leq 4\epsilon_{n-1}^{\frac{1}{10}}$. Thus $(\ref{estisharp})$ is fulfilled.

If it is in the non-resonant case, assume $j_0$ is the last resonant step before $n$.

If $j_0$ exists, we have
$$
\| B_{j_0}\|_0\leq \epsilon_{j_0-1}^{-\frac{1}{800}},
$$
$$
A_{j_0}=e^{A_{j_0}''}, \ \ \lVert A_{j_0}''\rVert < 2\epsilon_{j_0-1}^{\frac{1}{10}}, \ \ \lVert A_{j_0}\rVert \leq 1+4\epsilon_{j_0-1}^{\frac{1}{10}}.
$$
By our choice of $j_0$, from $j_0$ to $n$, every step is non-resonant. Thus we have
\begin{equation}
\label{estnn}\lVert A_{n}- A_{j_0}\rVert\leq 4\epsilon_{j_0}^{\frac{1}{2}},
\end{equation}
so
$$
\lVert A_{n}\rVert \leq 1+4\epsilon_{j_0-1}^{\frac{1}{10}}+ 4\epsilon_{j_0}^{\frac{1}{2}}\leq 2\lVert A_0\rVert.
$$
Estimate $(\ref{estnn})$ implies that if we rewrite $A_{n}=e^{A_{n}''}$, then
$$
\lVert A_{n}''\rVert \leq 4\epsilon_{j_0-1}^{\frac{1}{10}}.
$$
Moreover, we have
$$
\lVert B_{n}\rVert_0 \leq 2\|B_{j_0}\|_0\leq 2\epsilon_{j_0-1}^{-\frac{1}{800}}.
$$
Similarly to the process of $(\ref{estt})$, $(\ref{estisharp})$ is fulfilled.

If $j$ vanishes, it immediately implies that from $1$ to $n$, each step is non-resonant. In this case, $\lVert A_{n}\rVert \leq 2\lVert A_0\rVert$ and the estimate $(\ref{estisharp})$ is naturally satisfied as
$$\lVert B_n\rVert_{0}\leq 2.$$
Thus, we finish the proof.
\end{pf}

\section{Proof of Theorem \ref{main}}
We denote by
$$
P_k(E)=\sum\limits_{j=1}^k (({S_E^{\lambda V}(\theta)})_{2j-1})^*({S_E^{\lambda V}(\theta)})_{2j-1}.
$$
To prove Theorem \ref{main}, we only need to prove the following Lemma.
\begin{Lemma}\label{P}
Assume $\alpha\in DC_d$ and $V\in C^\omega(\T^d,\R)$, there exists $\lambda_0(\alpha,V)$ and $C=C(\alpha,\lambda V)$ such that if $\lambda<\lambda_0$, then for any $E\in \Sigma_{\alpha,\lambda V}$, we have $\|P_k(E)\|\leq C\|P^{-1}_k(E)\|^{-3}$.
\end{Lemma}
\noindent\textbf{Proof of Theorem \ref{main}:} It is standard that $\|P_k(E)\|=\det{P_k(E)}\|P^{-1}_k(E)\|$ since $P_k(E)$ is a self-adjoint matrix. If we let $\det{P_k(E)}=\frac{1}{4\epsilon_k^2}$, by Lemma \ref{P}, we have
$$
\|P_k(E)\|=\frac{1}{4\epsilon_k^2}\|P^{-1}_k(E)\|\leq \frac{C}{\epsilon_k^2}\|P_k(E)\|^{-\frac{1}{3}},
$$
thus
$$
\|P_k(E)\|\leq C\epsilon_k^{-\frac{3}{2}}.
$$
By Proposition \ref{mfunction}, we have
$$
\frac{\Im M(E+i\epsilon_k)}{\epsilon_k}\leq C\epsilon_k^{-\frac{3}{2}}.
$$

On the one hand, for any bounded potential $V$ and any solution $u$ we have $\|u\|_{2(k+1)}\leq C\|u\|_{2k}$, by Proposition \ref{prob} we have $\det{P_{k+1}(E)}\leq C\det{P_{k}(E)}$, thus $\epsilon_{k+1}\geq c\epsilon_k$. On the other hand, we can check easily in \eqref{Mfunction} that $\frac{\Im M(E+i\epsilon)}{\epsilon}$ is monotonic with respect to $\epsilon$, thus for any $\epsilon>0$, there exists $k$ such that $\epsilon_{k+1}<\epsilon<\epsilon_k$, combining this with the fact $\epsilon_{k+1}\geq c\epsilon_k$, we have
$$
\frac{\Im M(E+i\epsilon)}{\epsilon}\leq \frac{\Im M(E+i\epsilon_{k+1})}{\epsilon_{k+1}}\leq C\epsilon_{k+1}^{-\frac{3}{2}}\leq C\epsilon_k^{-\frac{3}{2}}\leq C\epsilon^{-\frac{3}{2}}.
$$

By \eqref{ine1}, we have
\begin{equation}\label{spc}
\mu_{\theta}(E-\epsilon,E+\epsilon)\leq 2\epsilon\Im M(E+\epsilon)\leq C\epsilon^{\frac{1}{2}},
\end{equation}
for $E\in \Sigma_{\alpha,\lambda V}$ and $\theta\in\T^d$.

For any $E\in \R$, we have the following two cases\\
\noindent \textbf{Case 1:} $(E-\epsilon,E+\epsilon)\cap\Sigma_{\alpha,\lambda V}=\emptyset$, we have
$$
\mu_{\theta}(E-\epsilon,E+\epsilon)=0\leq C\epsilon^{\frac{1}{2}}.
$$
\noindent \textbf{Case 2:} $(E-\epsilon,E+\epsilon)\cap\Sigma_{\alpha,\lambda V}\neq\emptyset$, there exists $E'\in (E-\epsilon,E+\epsilon)\cap\Sigma_{\alpha,\lambda V}$, then
$$
(E-\epsilon,E+\epsilon)\subset(E'-2\epsilon,E'+2\epsilon).
$$
Thus
$$
\mu_{\theta}(E-\epsilon,E+\epsilon)\leq\mu_{\theta}(E'-2\epsilon,E'+2\epsilon)\leq C(2\epsilon)^{\frac{1}{2}}\leq C\epsilon^{\frac{1}{2}}.
$$
\noindent\textbf{Proof of Lemma \ref{P}:}
Note that
$$
S_E^{\lambda V}(\theta)=\begin{pmatrix}E&-1\\1&0\end{pmatrix}(I+\begin{pmatrix}0&0\\ \lambda V(\theta)&0\end{pmatrix}),
$$
thus there exist $\lambda_0(\alpha,V)>0$ and $f_0\in C^\omega_{r_0}(\T^d,sl(2,\R))$ such that if $\lambda<\lambda_0$, we can rewrite
$$
I+\begin{pmatrix}0&0\\ \lambda V(\theta)&0\end{pmatrix}=e^{f_0(\theta)},
$$
with $\|f_0\|_{r_0}\leq \epsilon_*$ where $\epsilon_*$ is defined in Proposition \ref{reducibility}.

Hence we have $S_E^{\lambda V}(\theta)=A_0e^{f_0(\theta)}$ with $A_0=\begin{pmatrix}E&-1\\1&0\end{pmatrix}$. Since $\alpha\in DC_d$ and $V\in C^\omega(\T^d,\R)$, by Proposition \ref{reducibility}, there exist $B_n\in C_{r}^\omega(\T^d, PSL(2,\R))$ and $A_n\in SL(2,\R)$ satisfying
$$
B_n^{-1}(\theta+\alpha)A_0e^{f_0(\theta)}B_n(\theta)=A_ne^{f_n(\theta)},
$$
with estimates
\begin{equation*}
\|f_n\|_{r}\leq \epsilon_n,\ \ \|B_n\|_{0}\leq \epsilon_{n-1}^{-\frac{1}{800}}.
\end{equation*}
Moreover,
there exists unitary $U_n\in SL(2,\C)$ such that
$$
U_nA_{n}U_n^{-1}=\begin{pmatrix} e^{i\xi_n} & c_n\\ 0 & e^{-i\xi_n} \end{pmatrix},
$$
and
\begin{equation}
\lvert c_n\rvert \lVert B_{n}(\theta)\rVert_0^{8}\leq 4\lVert A_0\rVert,
\end{equation}
with $\xi_n,c_n\in \C$.

For $E\in\Sigma_{\alpha,\lambda V}$, we always have that $|\Im\xi_n|\leq \epsilon_n^{\frac{1}{4}}$ since $(\alpha,A_ne^{f_n})$ is not uniformly hyperbolic. Let $\Phi_n(\theta)=B_n(\theta)U_n^{-1}\in C^\omega_r(\T^d,PSL(2,\C))$ we have
\begin{align}\label{con}
\Phi_n^{-1}(\theta+\alpha)S_E^{\lambda V}(\theta)\Phi_n(\theta)=\tilde{A}_ne^{\tilde{f}_n(\theta)},
\end{align}
where $\tilde{A}_n=\begin{pmatrix}e^{2\pi i\gamma_n}&\tilde{c}_n\\0&e^{-2\pi i\gamma_n}\end{pmatrix}$ with $\gamma_n=\frac{\Re\xi_n}{2\pi}$ and
\begin{equation}\label{esta1}
\|\tilde{f}_n\|_r\leq\epsilon_n^{\frac{1}{4}},
\end{equation}
\begin{equation}\label{esta3}
\|\Phi_n\|_0\leq \|B_n\|_0\|U_n^{-1}\|=\|B_n\|_0\leq\epsilon_{n-1}^{-\frac{1}{800}},
\end{equation}
\begin{equation}\label{esta2}
|\tilde{c}_n|\|\Phi_n\|^8_0\leq 4\|A_0\| \leq C.
\end{equation}

For any $k\in (\epsilon_0^{-\frac{1}{40}},\infty)$, we denote by
$$
X_k=\sum\limits_{j=1}^k (\tilde{A}^{2j-1}_n)^*\tilde{A}^{2j-1}_n, \ \ k\in I_n:=(\epsilon_{n-1}^{-\frac{1}{40}},\epsilon_n^{-\frac{1}{30}}),
$$
$$
\tilde{X}_k(\theta)=\sum\limits_{j=1}^k((\tilde{A}_ne^{\tilde{f}_n(\theta)})_{2j-1})^*(\tilde{A}_ne^{\tilde{f}_n(\theta)})_{2j-1}, \ \ k\in I_n:=(\epsilon_{n-1}^{-\frac{1}{40}},\epsilon_n^{-\frac{1}{30}}).
$$

We divide the remaining proof into three steps.\\
\textbf{STEP 1: Estimation of $X_k$.}

By Proposition \ref{comp}, we have
$$
X_k=\begin{pmatrix}
k&x_{k,1}\\
\bar{x}_{k,1}&x_{k,2}
\end{pmatrix}
$$
where
\begin{align}\label{e1}
x_{k,1}=\tilde{c}_ne^{-2\pi i\gamma_n}\sum\limits_{j=1}^k\frac{e^{-4\pi i\gamma_n(2j-1)}-1}{e^{-4\pi i \gamma_n}-1},
\end{align}
\begin{align}\label{e2}
x_{k,2}=k+|\tilde{c}_n|^2\sum\limits_{j=1}^k(\frac{\sin2\pi(2j-1)\gamma_n}{\sin2\pi\gamma_n})^2.
\end{align}
By \eqref{e1} and \eqref{e2}, we have
\begin{equation}\label{eq1}
x_{k,1}=\frac{\tilde{c}_ne^{-2\pi i\gamma_n}}{e^{-4\pi i\gamma_n}-1}(e^{-4\pi i\gamma_n}\frac{e^{-8\pi ik\gamma_n}-1}{e^{-8\pi i\gamma_n}-1}-k)
\end{equation}
\begin{equation}\label{es2}
x_{k,2}=k(1+\frac{2|\tilde{c}_n|^2}{|e^{-4\pi i\gamma_n}-1|^2}(1-\frac{\sin8\pi k\gamma_n}{2k\sin4\pi \gamma_n})),
\end{equation}
\begin{equation}\label{es3}
\det{X_k}=k^2(1+\frac{|\tilde{c}_n|^2}{|e^{-4\pi i\gamma_n}-1|^2}(1-(\frac{\sin4\pi k\gamma_n}{k\sin4\pi\gamma_n})^2)).
\end{equation}

Note that by \eqref{es2}, we have $x_{k,2}>k$. Since $X_{k}$ is positive we have
$$
kx_{k,2}>(x_{k,1})^2,
$$
this implies that
$$
x_{k,2}>x_{k,1}.
$$
It is easy to see that
$$
x_{k,2}\leq\|X_k\|\leq 2x_{k,2}.
$$
Thus
\begin{align*}
\|(X_k)^{-1}\|^{-1}=\frac{\det{X_k}}{\|X_k\|}\geq\frac{\det{X_k}}{2x_{k,2}}.
\end{align*}

By \eqref{es2} and \eqref{es3}, we have
\begin{align*}
\|(X_k)^{-1}\|^{-1}&\geq\frac{k(1+\frac{|\tilde{c}_n|^2}{|e^{-4\pi i\gamma_n}-1|^2}(1-(\frac{\sin4\pi k\gamma_n}{k\sin4\pi\gamma_n})^2))}{2(1+\frac{2|\tilde{c}_n|^2}{|e^{-4\pi i\gamma_n}-1|^2}(1-\frac{\sin8\pi k\gamma_n}{2k\sin4\pi \gamma_n}))}.
\end{align*}
We have the following two cases:\\
Case 1: $k\|4\gamma_n\|_{\R/\Z}\geq\frac{2}{3}$, we have
$$
1-(\frac{\sin4\pi k\gamma_n}{k\sin4\pi\gamma_n})^2\geq \frac{1}{4},
$$
$$
1-\frac{\sin8\pi k\gamma_n}{2k\sin4\pi \gamma_n}\leq 2,
$$
thus
$$
k(1+\frac{|\tilde{c}_n|^2}{|e^{-4\pi i\gamma_n}-1|^2}(1-(\frac{\sin4\pi k\gamma_n}{k\sin4\pi\gamma_n})^2))\geq \frac{k}{4}(1+\frac{|\tilde{c}_n|^2}{|e^{-4\pi i\gamma_n}-1|^2}),
$$
$$
2(1+\frac{2|\tilde{c}_n|^2}{|e^{-4\pi i\gamma_n}-1|^2}(1-\frac{\sin8\pi k\gamma_n}{2k\sin4\pi \gamma_n}))\leq 8(1+\frac{|\tilde{c}_n|^2}{|e^{-4\pi i\gamma_n}-1|^2}).
$$
this implies
$$
\|(X^n_k)^{-1}\|^{-1}\geq c k.
$$
\noindent Case 2: $k\|4\gamma_n\|_{\R/\Z}\leq\frac{2}{3}$, we have
$$
1-(\frac{\sin4\pi k\gamma_n}{k\sin4\pi\gamma_n})^2\geq \frac{1}{4}k^2\|4\gamma_n\|^2_{\R/\Z},
$$
$$
1-\frac{\sin8\pi k\gamma_n}{2k\sin4\pi \gamma_n}\leq 100k^2\|4\gamma_n\|^2_{\R/\Z},
$$
thus
$$
k(1+\frac{|\tilde{c}_n|^2}{|e^{-4\pi i\gamma_n}-1|^2}(1-(\frac{\sin4\pi k\gamma_n}{k\sin4\pi\gamma_n})^2))\geq \frac{k}{4}(1+\frac{|\tilde{c}_n|^2}{|e^{-4\pi i\gamma_n}-1|^2}k^2\|4\gamma_n\|^2_{\R/\Z}),
$$
$$
2(1+\frac{2|\tilde{c}_n|^2}{|e^{-4\pi i\gamma_n}-1|^2}(1-\frac{\sin8\pi k\gamma_n}{2k\sin4\pi \gamma_n}))\leq 400(1+\frac{|\tilde{c}_n|^2}{|e^{-4\pi i\gamma_n}-1|^2}k^2\|4\gamma_n\|^2_{\R/\Z}).
$$
this implies
\begin{equation*}
\|(X_k)^{-1}\|^{-1}\geq c k.
\end{equation*}

Thus for any $k\in I_n$, we always have
\begin{equation}\label{eq1}
\|(X_k)^{-1}\|^{-1}\geq c k,
\end{equation}
\begin{equation}\label{eq2}
\|X_k\|\leq 2x_{k,2}\leq Ck(1+k^2|\tilde{c}_n|^2).
\end{equation}
\textbf{STEP 2: Estimation of $\tilde{X}_k(\theta)$.}

We need the following Lemma in \cite{aj},
\begin{Lemma}[\cite{aj}]\label{aj}
Assume $T=\begin{pmatrix}e^{2\pi i\theta}&c\\0&e^{-2\pi i\theta}\end{pmatrix}\in SL(2,\C)$ and $\tilde{T}\in C^0(\T^d,SL(2,\C))$, let $\tilde{T}_k(x)=\sum\limits_{j=1}^k\tilde{T}_{2j-1}^*(x)\tilde{T}_{2j-1}(x)$ and $T_k=\sum\limits_{j=1}^k (T^{2j-1})^*T^{2j-1}$, if $\|\tilde{T}-T\|_0\leq \frac{1}{100}k^{-2}(1+2ck)^{-2}$, we have
$$
\|\tilde{T}_k-T_k\|_0\leq 1.
$$
\end{Lemma}
Note that by \eqref{esta1} we have $\|\tilde{A}_ne^{\tilde{f}_n}-\tilde{A}_n\|_0\leq 2\epsilon_n^{\frac{1}{4}}$. Since $k\in I_n$, we have $\epsilon_n^{\frac{1}{4}}\leq \frac{1}{100}k^{-2}(1+2\tilde{c}_nk)^{-2}$, by Lemma \ref{aj}, we have
$$
\|\tilde{X}_k-X_k\|_0\leq 1,
$$
thus
$$
\|(\tilde{X}_k)^{-1}-(X_k)^{-1}\|_0\leq \|(\tilde{X}_k)^{-1}\|_0\|\tilde{X}_k-X_k\|_0\|(X_k)^{-1}\|_0\leq 1.
$$
By \eqref{eq1} and \eqref{eq2}, for any $k\in I_n$, we have
$$
\|(\tilde{X}_k)^{-1}\|^{-1}\geq c k.
$$
$$
\|\tilde{X}_k\|\leq Ck(1+k^2|\tilde{c}_n|^2).
$$
\textbf{STEP 3: Estimation of $P_k(E)$.}

For $k\in I_n$, by equation \eqref{con}, we have
\begin{align*}
\|P_k(E)\|_0\leq \|\Phi_n\|_0^4\|\tilde{X}_k\|_0&\leq C\|\Phi_n\|_0^4k(1+k^2|\tilde{c}_n|^2),
\end{align*}
\begin{align*}
\|P^{-1}_k(E)\|_0^{-1}\geq \|\Phi_n\|_0^{-4}\|(\tilde{X}_k)^{-1}\|_0^{-1}&\geq c\|\Phi_n\|_0^{-4}k,
\end{align*}
thus
$$
\frac{\|P_k(E)\|_0}{\|P^{-1}_k(E)\|_0^{-3}}\leq C\|\Phi_n\|_0^{16}|\tilde{c}_n|^2+C\|\Phi_n\|_0^{16}k^{-2}.
$$

On the one hand, $k^{-2}\leq \epsilon_{n-1}^{\frac{1}{20}}$, by \eqref{esta3}, we have $C\|\Phi_n\|_0^{16}k^{-2}\leq C$. On the other hand, by \eqref{esta2}, we have $C\|\Phi_n\|_0^{16}|\tilde{c}_n|^2\leq C$, thus for any $k\in (\epsilon_0^{-\frac{1}{40}},\infty)$
$$
\frac{\|P_k(E)\|_0}{\|P^{-1}_k(E)\|_0^{-3}}\leq C.
$$

For $k\in (0,\epsilon_0^{-\frac{1}{40}})$, it is obvious that there exists $C=C(\alpha,\lambda V)$ such that
$$
\frac{\|P_k(E)\|_0}{\|P^{-1}_k(E)\|_0^{-3}}\leq C.
$$

Thus we finish the proof of Lemma \ref{P}. \qed
\section{Appendix}
In this section, we give an iteration proposition proved in \cite{lyzz,ccyz} as a generalization of the results in \cite{hy}.
\begin{Proposition}\label{prop1}
Let $\alpha\in DC_d(\kappa,\tau)$, $\kappa$, $r>0$, $\tau>d-1$.
Suppose that $A\in SL(2,\R)$, $f\in C^{\omega}_{r}(\T^{d},sl(2,\R))$.  Then for any $r'\in (0,r)$, there exist $c=c(\kappa,\tau,d)$ and a numerical constant $D$ such that if
\begin{equation}\label{estf}
\lvert f \rvert_r\leq\epsilon \leq \frac{c}{\lVert A\rVert^D}(r-r')^{D\tau},
\end{equation}
then there exist $B\in C^{\omega}_{r'}(\T^{d},PSL(2,\R))$, $A_{+}\in SL(2,\R)$ and $f_{+}\in C^{\omega}_{r'}(\T^{d},$
$sl(2,\R))$ such that
$$
B^{-1}(\theta+\alpha)(Ae^{f(\theta)})B(\theta)=A_{+}e^{f_+(\theta)}.
$$
More precisely, let $N=\frac{2}{r-r'} \lvert \ln \epsilon \rvert$, then we can distinguish two cases:
\begin{itemize}
\item (Non-resonant case)   if for any $n\in \Z^{d}$ with $0<\lvert n \rvert \leq N$, we have
$$
\lvert 2\xi - <n,\alpha> \rvert\geq \epsilon^{\frac{1}{10}},
$$
then
$$\lvert B-Id\rvert_{r'}\leq \epsilon^{\frac{1}{2}} ,\   \ \lvert f_{+}\rvert_{r'}\leq 4\epsilon^{3-\frac{1}{5}}.$$
and
$$\lVert A_+-A\rVert\leq 2\lVert A\rVert\epsilon.$$
\item (Resonant case) if there exists $n_\ast$ with $0<\lvert n_\ast\rvert \leq N$ such that
$$
\lvert 2\xi- <n_\ast,\alpha> \rvert< \epsilon^{\frac{1}{10}},
$$
then
$$\lvert B \rvert_{r'}\leq \epsilon^{-\frac{1}{1600}}\times\epsilon^{\frac{-r'}{r-r'}},\ \ \lVert B\rVert_0 \leq \epsilon^{-\frac{1}{1600}},\ \ \lvert f_{+}\rvert_{r'}\ll \epsilon^{1600}.$$ Moreover, $A_+=e^{A''}$ with $\lVert A''\rVert \leq 2\epsilon^{\frac{1}{10}}$.
\end{itemize}
\end{Proposition}

%
%

\section*{Acknowledgement}
We are indebted to Professor Jiangong You and Professor Yiqian Wang for their enthusiastic help. We are
also grateful to Lingrui Ge for carefully reading the manuscript and many useful suggestions which gave us positive hints. X. Zhao was supported by China Scholarship Council (CSC)(No. 201906190072) and NSFC grant (11771205).

\end{document}